\documentclass[12pt,leqno]{article}

\usepackage{amsmath,amssymb,amsthm}

\makeatletter

\newtheorem{Theorem}{Theorem}
\newtheorem{Proposition}[Theorem]{Proposition}
\newtheorem{Lemma}[Theorem]{Lemma}
\newtheorem{Corollary}[Theorem]{Corollary}

\theoremstyle{definition}

\newtheorem{Def}[Theorem]{Definition}

\theoremstyle{remark}

\def\({{\rm (}}
\def\){{\rm )}}

\let\Mathrm\operator@font
\let\Cal\mathcal

\def\standop#1{\mathop{\Mathrm #1}\nolimits}
\def\difstop#1#2{\expandafter\def\csname #1\endcsname{\standop{#2}}}
\def\defstop#1{\difstop{#1}{#1}}

\defstop{codim}
\defstop{id}
\defstop{Sym}
\defstop{Hom}
\defstop{End}
\defstop{rank}
\defstop{Spec}
\def\GL{{\sl{GL}}}
\defstop{Ker}
\defstop{Min}
\difstop{Image}{Im}
\difstop{tdeg}{trans.deg}

\defstop{length}
\defstop{supp}
\defstop{Proj}
\defstop{Flat}
\defstop{Pfaff}
\defstop{Pf}
\defstop{Alt}
\defstop{Mat}
\difstop{height}{ht}

\def\specialarrow#1{\setbox\z@=\hbox{$\m@th
 \mathop{\vphantom{\rightarrow}}\limits^{\hspace{.5ex}{#1}\hspace
{.8ex}}$}\mathrel{\ifdim\wd\z@<1.2em\dimen\tw@
1.2em\else\dimen\tw@\wd\z@\fi\copy\z@\kern-\wd\z@\hbox to\dimen\tw@
{\rightarrowfill}}}

\def\sdarrow#1{\downarrow\hbox to 0pt{\scriptsize$#1$\hss}}
\def\suarrow#1{\uparrow\hbox to 0pt{\scriptsize$#1$\hss}}
\def\ssearrow#1{\searrow\hbox to 0pt{\scriptsize$#1$\hss}}

\def\ext{{\textstyle\bigwedge}}

\def\section{\@startsection{section}{1}{\z@ }%
{-3.5ex plus -1ex minus -.2ex}{2.3ex plus .2ex}{\bf }}

\long\def\refname{\par\kern -3ex
\begin{center}\rm R\sc{eferences}\end{center}\par\kern 
-2ex}

\def\@seccntformat#1{\csname the#1\endcsname.\quad}

\def\@@@sect#1#2#3#4#5#6[#7]#8{%
   \ifnum #2>\c@secnumdepth 
      \def \@svsec {}\else \refstepcounter {#1}%
      \def\@svsec{}
   \fi 
   \@tempskipa #5\relax 
   \ifdim \@tempskipa >\z@ 
     \begingroup #6\relax \@hangfrom {\hskip #3\relax 
     \@svsec}{\interlinepenalty \@M #8\par }\endgroup 
     \csname #1mark\endcsname {#7}
   \else 
   \def \@svsechd {#6\hskip #3\@svsec #8\csname #1mark\endcsname {#7}}
   \fi \@xsect {#5}}

\def\@@@startsection#1#2#3#4#5#6{%
 \if@noskipsec \leavevmode \fi \par \@tempskipa #4\relax \@afterindenttrue 
 \ifdim \@tempskipa <\z@ \@tempskipa -\@tempskipa \@afterindentfalse 
 \fi \if@nobreak \everypar {}\else \addpenalty {\@secpenalty }\addvspace 
  {\@tempskipa }\fi \@ifstar {\@ssect {#3}{#4}{#5}{#6}}{\@dblarg 
  {\@@@sect {#1}{#2}{#3}{#4}{#5}{#6}}}}

\def\theparagraph{\thesection.\arabic{paragraph}}
\def\aparagraph{\@@@startsection{paragraph}{2}{\z@ }%
              {1.75ex plus .2ex minus .15ex}{-1em}{\bf(\theparagraph) } }
\def\paragraph{\@@@startsection{paragraph}{2}{\z@ }%
              {1.75ex plus .2ex minus .15ex}{-1em}{}{\bf(\theparagraph)} }

\let\c@Theorem\c@paragraph

\title{Another proof of theorems of De Concini and Procesi\thanks{%
{\em 2000 Mathematics Subject Classifications.}
Primary 13A50; Secondary 14L30.\newline
{\em Key Words and Phrases.} invariant theory, determinantal ideal, 
Pfaffian ideal, Cohen--Macaulay}%
\\
\small Dedicated to Professor Kei-ichi Watanabe on the occasion of his
sixtieth birthday} 

\author{M{\sc itsuyasu} H{\sc ashimoto}}
\date{\normalsize
Graduate School of Mathematics, Nagoya University\\
Chikusa-ku,  Nagoya 464--8602 JAPAN\\
\small \tt hasimoto@math.nagoya-u.ac.jp}

\begin{document}

\maketitle

\begin{abstract}
We give a new proof of 
some characteristic-free fundamental theorems in invariant theory first
proved in 
C. De Concini and C. Procesi, A characteristic free approach
to invariant theory, {\em Adv. Math.} {\bf 21} (1976), 330--354.
We treat the action of the general linear group and the symplectic group.
Our approach is geometric, and 
utilizes the fact that the categorical quotients are
principal fiber bundles off codimension two or more.
\end{abstract}

\section{Introduction}

Let $k$ be an algebraically closed field, and $n\geq m\geq t\geq 2$.
Set $r=t-1$, $E=k^r$, $V=k^n$, and $W=k^m$.
Let $M=\Hom(E,W)\times\Hom(V,E)$, and $G=\GL(E)$.
$G$ acts on $M$ by $g(\varphi,\psi)=(\varphi g^{-1},g\psi)$.
Let $Y_t:=\{f\in\Hom(V,W)\mid \rank f<t\}$.
It is easy to see that $\pi\colon M\rightarrow Y_t$ given by
$\pi(\varphi,\psi)=\varphi\circ \psi$ is well-defined and $G$-invariant.

De Concini and Procesi \cite{DP} proved that $\pi$ is a categorical quotient.
The case of characteristic zero was proved in \cite{Weyl}.
In other words, $\pi$ induces an isomorphism $k[Y_t]\cong k[M]^G$.
Yet another interpretation is as follows.
$M$ is isomorphic to $E^{\oplus n}\oplus (E^*)^{\oplus m}$ as a $G$-module,
where $G$ acts via
\[
g(e_1,\dots,e_n,e_1^*,\dots,e_m^*)=
(ge_1,\dots,ge_n,e_1^*g^{-1},\dots,e_m^*g^{-1}).
\]
It is easy to see that $\xi_{ij}=e^*_i(e_j)$ is a $G$-invariant polynomial
function, i.e., an element in $k[M]^G$.
The De Concini--Procesi theorem says that $k[M]^G$ is generated by $\xi_{ij}$,
and the kernel of the surjective map $k[x_{ij}]\rightarrow k[M]^G$ given by
$x_{ij}\mapsto \xi_{ij}$ is generated by the determinantal ideal $I_t(x_{ij})$
generated by the all $t$-minors of the matrix $(x_{ij})$.

The purpose of this article is to give a short and geometric proof to
the theorem.
Here we give a sketch of the proof.

We need to assume that $Y_t$ is Cohen--Macaulay.
This was first proved by Hochster and Eagon \cite{HE}.
See also \cite{BV} and \cite{BC}.
The Cohen--Macaulay property of $Y_t$ is highly non-trivial, so our
proof is not completely self-contained.
But the rest of the argument is easy and self-contained.

It is easy to see that $Y_t$ satisfies the $(R_1)$ condition, so 
$Y_t$ is normal.
As the codimension of $Y_{t-1}$ in $Y_t$ is at least two, when we set
$U=Y_{t}\setminus Y_{t-1}$, we have $k[Y_t]=\Gamma(U,\Cal O_U)$.
On the other hand, $\pi|_{\pi^{-1}(U)}\colon \pi^{-1}(U)\rightarrow U$ is
a principal $G$-bundle.
Hence $\Gamma(U,\Cal O_U)\rightarrow \Gamma(\pi^{-1}(U),\Cal O_M)^G$ is
an isomorphism.
This is enough to prove the theorem.

In section~2, we prepare necessary results to prove our main theorem.
In section~3, we prove the main theorem as described above.
In section~4, we extend the result to arbitrary base ring, rather than
an algebraically closed base field.
The key is good filtrations of representations of algebraic groups.
In section~5, we show that a similar argument is effective for
the symplectic group action.

The results proved in section~4 and 5 are also already proved in \cite{DP},
and there is no new result in this paper.
What is new here is a simple geometric proof of their results 
based on the Cohen--Macaulay
property of the candidates of the invariant subrings.

\section{Preliminaries}\label{preliminaries}

Let $k$ be an algebraically closed field.

In the sequel, we will apply set theoretic argument to $k$-varieties.
For example, given a $k$-morphism $f:X\rightarrow Y$, we will say that 
$f$ is an isomorphism, since 
the induced map of the $k$-valued points $f(k):X(k)\rightarrow Y(k)$ is
bijective.
This argument itself is obviously false (for example, consider 
$\Spec k[t]\rightarrow \Spec k[t^2,t^3]$),
but the reader should interpret it into
the correct argument: as for any finitely generated $k$-algebra $A$, the
induced map of the set of $A$-valued points $f(A):X(A)\rightarrow Y(A)$
is bijective, $f$ is an isomorphism.
We will abuse such a ``$k$-valued points only'' argument for brevity,
only in the
case where the interpretation to the correct argument is straightforward
(but annoying).

Let $m,n,t$ be integers such that $2\leq t\leq m\leq n$.
Set $V:=k^n$, $W:=k^m$, and $E:=k^r$, where $r:=t-1$.
Set $M:=\Hom(E,W)\times \Hom(V,E)$, and 
$Y_t:=Y_t(V,W)=\{f\in\Hom(V,W)\mid \rank f < t\}\subset \Hom(V,W)$.
Let $S=k[x_{ij}]_{1\leq i\leq m,\,1\leq j\leq n}$ be the polynomial ring in
$mn$ variables over $k$, so that $S$ is the coordinate ring of $\Hom(V,W)$.
Then $Y_t$ is a closed subscheme of $\Hom(V,W)$ defined by $I_t$, 
where $I_t=I_t(x_{ij})$ is the ideal of $S$
generated by $t$-minors of the matrix $(x_{ij})$.
Let $T:=k[w_{il},v_{lj}\mid 1\leq i\leq m,\;1\leq l<t,\;1\leq j\leq n]$.
Then $T$ is the coordinate ring of $M$.

We define $\pi: M\rightarrow Y_t$ by $\pi(\varphi,\psi)=\varphi\circ \psi$
for $\varphi\in\Hom(E,W)$ and $\psi\in\Hom(V,E)$.
This is a well-defined morphism, since a linear map which factors through
the $r$-dimensional space $E$ has rank less than $t$.
The associated map of the coordinate rings is given by $x_{ij}\mapsto
\sum_{l=1}^{t-1}w_{il}v_{lj}$.

The following argument is taken from \cite[pp.4--5]{BV} for convenience of
readers.

Obviously, $\pi$ is surjective.
As $M$ is irreducible, $Y_t$ is also irreducible.

Let $V=U\oplus \tilde U$ be a direct sum decomposition with $\dim U=r$.
If $f\in Y_t$ and $f|_U$ is injective, then there exist unique linear
maps $g\colon \tilde U\rightarrow U$ and $h\colon U\rightarrow W$ such that
$f(u\oplus \tilde u)=h(u)+h(g(\tilde u))$ for all $u\in U$ and $\tilde u
\in \tilde U$.
So if we set 
\[
N:=\{f\in Y_t\mid f|_U \text{ injective}\},
\]
then there is an isomorphism of $k$-schemes
\[
\Hom(\tilde U,U)\times (\Hom(U,W)\setminus Y_{r-1}(U,W))\rightarrow N.
\]
Since the variety on the left is an open subvariety of 
$\Hom(\tilde U,U)\times \Hom(U,W)$, we have that
\begin{multline*}
\dim Y_t(V,W)=\dim N=\dim \Hom(\tilde U,U)\times \Hom(U,W)\\
=(m-r)r+rn=mr+nr-r^2.
\end{multline*}
Moreover, $N$ is non-singular.
Varying $U$, we have that $Y_t(V,W)\setminus Y_{t-1}(V,W)$ is non-singular.
So we have the following.

\def\citinfo{\cite[(1.1)]{BV}}
\begin{Proposition}[\citinfo]
\begin{description}
\item[1] $Y_t(V,W)$ is irreducible.
\item[2] $\dim Y_t(V,W)= mr+nr-r^2$.
\item[3] $Y_t$ is non-singular off $Y_{t-1}$.
\end{description}
\end{Proposition}

Hence we have

\begin{Lemma}\label{codim3.thm}
$\dim Y_t-\dim Y_{t-1}=m+n-2r+1\geq 3$.
In particular, $Y_t$ satisfies the $(R_2)$ condition.
\end{Lemma}

We need the following theorem, which was first proved by Hochster and
Eagon \cite{HE}.
See also \cite{BV}.

\begin{Theorem}
$Y_t$ is Cohen--Macaulay.
\end{Theorem}

Since $Y_t$ is irreducible and satisfies the $(R_1)$ and the $(S_2)$ 
conditions, we have the following immediately.

\begin{Corollary}
$Y_t$ is normal and integral.
\end{Corollary}

\begin{Def}
Let $G$ be an affine algebraic group over $k$, $X$ a $G$-action which
is of finite type over $k$, and
$f\colon X\rightarrow Y$ a $k$-morphism.
We say that $f$ is $G$-invariant if $f(gx)=f(x)$ holds for $x\in X$ and 
$g\in G$.
We say that $f$ is a principal $G$-bundle if $f$ is faithfully flat, 
$G$-invariant, and the map $\Phi\colon G\times X \rightarrow X\times_Y X$
given by $\Phi(g,x)=(gx,x)$ for $g\in G$ and $x\in X$ is an isomorphism.
\end{Def}

It is not so difficult to show that a principal $G$-bundle is a 
universally submersive geometric quotient in the sense of \cite{GIT}.
We do not prove this because we will not use it later.
What we need is the following.

\begin{Lemma}\label{torsor.thm}
Let $G$ be an affine algebraic group over $k$, $X$ a $G$-action which
is of finite type over $k$, and
$\pi\colon X\rightarrow Y$ a principal $G$-bundle.
Then the canonical map $\Cal O_Y\rightarrow (\pi_*\Cal O_X)^G$ is an
isomorphism.
\end{Lemma}

\proof First, consider the case $X=G\times Y$ and $\pi$ is the 
second projection.
The question is local on $Y$, so we may assume that $Y=\Spec A$ is affine.
Then the assertion reads $A\rightarrow (k[G]\otimes A)^G$ is an isomorphism,
and this is trivial.

Applying this observation to the second projection $G\times X'\rightarrow X'$,
and considering the $X'$-isomorphism $\Phi\colon G\times X'\cong X\times_Y X'$
which is also a $G$-isomorphism,
we have $\Cal O_{X'}\rightarrow ((p_2)_*\Cal O_{X\times_Y X'})^G$ is an 
isomorphism, where $X'$ is the scheme $X$ with the trivial $G$-action,
and $p_2\colon X\times_Y X' \rightarrow X'$ is the second projection.

Now apply $(\pi')^*$ to $\Cal O_Y\rightarrow (\pi_*\Cal O_X)^G$, where
$\pi'\colon X'\rightarrow Y$ is $\pi$ (remember that $X'=X$).
Since the $G$-invariance is compatible with the flat base change,
the result is
\[
\Cal O_{X'}\rightarrow ((\pi')^*\pi_*\Cal O_X)^G\cong 
((p_2)_*\Cal O_{X\times_Y X'})^G.
\]
We know that this is an isomorphism.
As $\pi'$ is faithfully flat, we have that 
$\Cal O_Y\rightarrow (\pi_*\Cal O_X)^G$ is also an isomorphism.
\qed

\section{Main Theorem}

Let $k$ be an algebraically closed field, and $m,n,t\in\Bbb Z$ such that
$2\leq t\leq m\leq n$.
Set $V:=k^n$, $W:=k^m$, and $E:=k^{r}$, where $r=t-1$.
As in section~\ref{preliminaries}, set 
$M:=\Hom(E,W)\times\Hom(V,E)$, and $Y_t:=Y_t(V,W)=\{f\in\Hom(V,W)\mid
\rank f<t\}$.
Consider the morphism $\pi\colon M\rightarrow Y_t$ given by 
$\pi(\varphi,\psi)=\varphi\circ \psi$.
Denote by associated $k$-algebra map $S/I_t\rightarrow T$ by $\phi^\#$,
where $S/I_t$ and $T$ are as in section~\ref{preliminaries}.

Let $G:=\GL(E)$.
Then $G$ acts on $M$ via $g\cdot(\varphi,\psi)=(\varphi g^{-1}, g\psi)$
for $g\in G$, $\varphi\in\Hom(E,W)$, and $\psi\in\Hom(V,E)$.
Then obviously, $\pi$ is $G$-invariant.

The objective of this section is to give a new proof to the following theorem.

\def\citinfo{De Concini--Procesi \cite[section~3]{DP}}
\begin{Theorem}[\citinfo]\label{dp.thm}
Let the notation be as above.
Then the associated map $\phi^\#$ is injective, and the image of $\phi^\#$ 
is identified with $T^G$.
\end{Theorem}

\proof Since $\pi\colon M\rightarrow Y_t$ is dominating and $Y_t$ is 
integral, $\phi^\#$ is injective.
As $\pi$ is $G$-invariant, an injective map $\phi^\#\colon k[Y_t]=S/I_t
\rightarrow k[M]^G$ is induced.
It suffices to prove that $\phi^\#$ is surjective.

Set $U:=Y_t\setminus Y_{t-1}$.
Then since $Y_t$ is normal, 
we have that $\Gamma(U,\Cal O_U)=k[Y_t]$
by Lemma~\ref{codim3.thm}.

On the other hand, we claim that $\Gamma(U,\Cal O_U)=
\Gamma(\pi^{-1}(U),\Cal O_M)^G$.
If the claim is true, $\phi^\#$ is surjective and the proof is complete,
since $\Gamma(\pi^{-1}(U),\Cal O_M)^G\supset k[M]^G$.
Set $\tilde U=\pi^{-1}(U)$.
To prove the claim, it suffices to show that
\[
\pi|_{\tilde U}\colon \tilde U\rightarrow U
\]
is a principal 
$G$-bundle by Lemma~\ref{torsor.thm}.

Note that $\tilde U=\{(\varphi,\psi)\in M\mid \psi \text{ surjective and }
\varphi \text{ injective}\}$.

It is obvious that $\pi|_{\tilde U}$ is $G$-invariant, since $\pi$ is.

We prove that $\Phi\colon G\times \tilde U\rightarrow \tilde U
\times_U \tilde U$ given by $\Phi(g,\tilde u)=(g\tilde u,\tilde u)$ 
is an isomorphism.
Let $((\varphi,\psi),(\varphi',\psi'))\in\tilde U\times_U \tilde U$.
Since $\varphi\psi=\varphi'\psi'$ and $\varphi$ and $\varphi'$ are injective,
we have that
\[
\Ker \psi=\Ker \varphi\psi=\Ker\varphi'\psi'=\Ker\psi'.
\]
By the homomorphism theorem, it is easy to see that there exists a 
unique $g\in G(k)$ such that $\psi=g\psi'$.
Since 
\[
\varphi'\psi'=\varphi\psi=\varphi g\psi'
\]
and $\psi'$ is surjective, we have that $\varphi'=\varphi g$.
So
\[
((\varphi,\psi),(\varphi',\psi'))=((\varphi'g^{-1},g\psi'),(\varphi',\psi'))
=\Phi(g,(\varphi',\psi')).
\]
Hence $\Phi$ is surjective.

Next, assume that $\Phi(g,(\varphi,\psi))=\Phi(g',(\varphi',\psi'))$ 
for $g,g'\in G$, $\varphi,\varphi'\in\Hom(E,W)$, and 
$\psi,\psi'\in\Hom(V,E)$.
Then obviously $(\varphi,\psi)=(\varphi',\psi')$.
Since $g\psi=g'\psi$ and $\psi$ is surjective, $g=g'$.
Hence $\Phi$ is injective.

Since $\Phi$ is bijective, $\Phi$ is an isomorphism.

Next, we prove that $\pi|_{\tilde U}$ is faithfully flat.
Since $\pi$ is surjective, $\pi|_{\tilde U}$ is surjective.
We only need to prove that $\pi|_{\tilde U}$ is flat.

Since $\Phi$ is an isomorphism, for each closed point $x\in\tilde U$, 
the morphism $G\rightarrow \pi^{-1}(\pi(x))$ given by $g\mapsto gx$ is
an isomorphism.
So each fiber is integral, and the dimension of fibers is constant.
Since both $U$ and $\tilde U$ are non-singular, the flatness of 
$\pi|_{\tilde U}$ now follows easily from 
\cite[Corollary to Theorem~23.1]{CRT}.
\qed

\section{Arbitrary base ring}\label{arbitrary.sec}

In this section, we extend Theorem~\ref{dp.thm} to an arbitrary base ring $R$,
rather than an algebraically closed base field $k$.

First consider the case where the base ring $R$ is an arbitrary field.
The map $S/I_t \rightarrow R[M]^{G_R}$ of graded $R$-algebras is an 
isomorphism if (and only if) it is an isomorphism after taking a
faithfully flat base change of the base ring $R$.
As a field extension is faithfully flat, the map in problem is an isomorphism
by Theorem~\ref{dp.thm}.

Next consider the base ring $R=\Bbb Z$.
Since $\Hom_{\Bbb Z}(E_{\Bbb Z},W_{\Bbb Z})\times 
\Hom_{\Bbb Z}(V_{\Bbb Z},E_{\Bbb Z})$ is isomorphic to
$E_{\Bbb Z}^{\oplus n}\oplus (E_{\Bbb Z}^*)^{\oplus m}$ 
as a $G_{\Bbb Z}$-module,
the coordinate ring
\[
\Bbb Z[M_{\Bbb Z}]\cong (\Sym E^*_{\Bbb Z})^{\otimes n}
\otimes (\Sym E_{\Bbb Z})^{\otimes m}
\]
has a good filtration as a $G_{\Bbb Z}$-module.
Indeed, it is well-known that $\Sym E$ and $\Sym E^*$ have good filtrations.
By Mathieu's tensor product theorem \cite{Mathieu} and 
\cite[Corollary~III.4.1.8]{Hashimoto}, the tensor product $\Bbb Z[M_{\Bbb Z}]$ 
is good.
In particular, $H^i(G_{\Bbb Z},\Bbb Z[M])=0$ for $i>0$.
By the universal coefficient theorem, the canonical map
$\Bbb Z[M]^{G_\Bbb Z}\otimes k\rightarrow k[M]^G$ 
is an isomorphism for any field $k$.
It follows that $S/I_t\otimes k\rightarrow \Bbb Z[M]^{G_{\Bbb Z}}\otimes k$ 
is an 
isomorphism for any $k$.
This shows that $S/I_t\rightarrow \Bbb Z[M]^{G_{\Bbb Z}}$ is an isomorphism,
since each homogeneous component of $S/I_t$ and $\Bbb Z[M]^{G_{\Bbb Z}}$ are
finitely generated $\Bbb Z$-modules, and a $\Bbb Z$-linear map between
finitely generated $\Bbb Z$-modules is an isomorphism if and only if its
base change to an arbitrary field is.

Now consider arbitrary base ring $R$.
By the universal coefficient theorem again, 
$\Bbb Z[M]^G \otimes R\rightarrow R[M]^{G_R}$ is an isomorphism.
Hence

\def\citinfo{De Concini--Procesi \cite[section~3]{DP}}
\begin{Theorem}[\citinfo]\label{dp2.thm}
Let $R$ be an arbitrary commutative ring.
The canonical map $S/I_t\rightarrow R[M]^{G_R}$ is an isomorphism.
\end{Theorem}

\section{Symplectic group action}

We apply the strategy above to another example.

Let $k$ be an algebraically closed field, 
$t$ and $n$ integers such that $4\leq 2t\leq n$, 
$V:=k^n$ and $E:=k^{2t-2}$.

Let $A=(a_{ij})$ be an alternating $2t\times 2t$ matrix over a commutative
ring $R$.
Namely, $a_{ji}=-a_{ij}$ and $a_{ii}=0$.
We define the Pfaffian of $A$ to be
\[
\Pfaff(A)=\sum_{\sigma\in\Gamma}(-1)^\sigma a_{\sigma 1\,\sigma 2}
a_{\sigma 3 \,\sigma 4}\cdots a_{\sigma(2t-1)\,\sigma(2t)}\in R,
\]
where
\[
\Gamma=\{\sigma\in\mathfrak S_{2t}\mid \sigma1<\sigma 3<\cdots<\sigma(2t-1),
\;\sigma (2i-1)<\sigma (2i)\;(1\leq i\leq t)\}.
\]
Note that $(\Pfaff(A))^2=\det A$ (on the other hand, if $A$ is an alternating
matrix of odd size, $\det A=0$).

Let $\langle\;,\;\rangle:E\times E\rightarrow k$ be a bilinear form
given by 
$(\langle e_i,e_j\rangle)=\tilde J$,
where 
\[
J=J_{t-1}=\begin{pmatrix}
       &   & &  & 1    \\
       &   & &\cdot & \\
       &   & \cdot & & \\
       &   \cdot & & & \\
     1 &  &   &  &
   \end{pmatrix}
\text{ and }
\tilde J=
\tilde J_{t-1}=
             \begin{pmatrix}
                   O & J \\
                  -J & O
                \end{pmatrix},
\]
and $e_1,\ldots,e_{2t-2}$ is the standard basis of $E$.
Note that the bilinear form $\langle\;,\;\rangle$ induces a $k$-linear
map $\rho:\ext^2 E\rightarrow k$ given by $\rho(v\wedge v')=\langle v,v'
\rangle$.
We define $G$ to be the symplectic group 
$\{\varphi\in\End(E)\mid \rho\circ\ext^2\varphi=\rho\}$.
Note that $\dim G=2(t-1)^2+t-1$, see \cite[p.~3]{Humphreys}.

Define $M$ to be the affine space $\Hom(V,E)$.
Note that $G$ acts on $M$ in a natural way.
Define $Y_t$ to be the variety of $n\times n$ alternating matrices such that
Pfaffians of all main $2t\times 2t$ submatrices vanish,
where a main $2t\times 2t$ submatrix of an alternating $n\times n$ 
matrix $A=(a_{ij})$ is an alternating $2t\times 2t$ matrix of the form
$(a_{i_u i_v})_{1\leq u,v\leq 2t}$, where $1\leq i_1<\cdots<i_{2t}\leq n$.

Thus $Y_t$ is a closed subscheme of the affine space 
$(\ext^2 V)^*$.
Define $\pi\colon M\rightarrow (\ext^2 V)^*$ by 
$\pi(\varphi)=\rho\circ \ext^2\varphi$.

Note that there is a well-defined Pfaffian map $\Pf_{2t}^V\colon 
\ext^{2t}V\rightarrow \Sym\ext^2 V=k[(\ext^2 V)^*]$ given by
$\Pf_{2t}^V(v_1\wedge\cdots \wedge v_{2t})=\Pfaff(v_i\wedge v_j)$.
Note that $A\in (\ext^2 V)^*$ lies in $Y_t$ if and only if the composite
\[
\ext^{2t}V\xrightarrow{\Pf_{2t}}\Sym \ext^2 V=k[(\ext^2 V)^*]\xrightarrow{A}k
\]
is zero.

Let $\varphi\in\Hom(V,E)$.
As the diagram
\[
\begin{array}{ccccccc}
 & & \ext^{2t} V & \xrightarrow{\Pf_{2t}^V} & \Sym \ext^2 V &
\xrightarrow{\pi(\varphi)} & k\\
 & & \sdarrow{\ext^{2t}\varphi} & & \sdarrow{\Sym\ext^2\varphi} & 
 & \sdarrow{\id_k}\\
0 &= & \ext^{2t} E &  \xrightarrow{\Pf_{2t}^E} & \Sym \ext^2 E &
\xrightarrow{\rho} & k
\end{array}
\]
commutes, $\pi$ factors through $Y_t$.
So we have the morphism $\pi\colon M\rightarrow Y_t$.
By the definition of the symplectic group $G$, $\pi$ is $G$-invariant.

Fix the standard basis $f_1,\ldots, f_n$ of $V$.
For $f\in (\ext^2 V)^*$, the alternating matrix $A(f)=(f(f_i\wedge f_j))$
corresponds.
We say that $\varphi\in \Hom(V,V^*)$ is alternating if the representation
matrix $(a_{ij})$ given by $\varphi(f_j)=\sum_i a_{ij}f_i^*$ is alternating,
where $f_1^*,\ldots,f_n^*$ is the dual basis of $f_1,\ldots,f_n$.
This notion is independent of the choice of basis of $V$.
Denote the space of alternating maps by $\Alt(V)\subset \Hom(V,V^*)$.
Corresponding to $f\in(\ext^2 V)^*$, an alternating map 
$a(f)\in\Alt(V)$ is given by $a(f)(v)=\sum_i f(f_i\wedge v)f_i^*$.
The representation matrix of $a(f)$ is $A(f)$.

Similarly, we fix the basis $e_1,\ldots, e_{2t-2}$ of $E$, and
corresponding to $h\in (\ext^2 E)^*$, the $(2t-2)\times (2t-2)$ 
alternating matrix $A(h)$ corresponds, and the alternating map $a(h)$
corresponds.
Note that $A(\rho)=\tilde J$.
We denote by $\tilde \rho\in\Alt(E)$ the corresponding map $a(\rho)$.

It is easy to see that $a(\pi(\varphi))=\varphi^*\circ \tilde\rho\circ\varphi$.
In other words, $A(\pi(\varphi))={}^t X \tilde J X$, where
$X$ is the representation matrix of $\varphi$.

\begin{Lemma}\label{alternating.thm}
Let $A$ be an $n\times n$ alternating matrix of rank $2r$.
Then there exists some $X\in\GL_n(k)$ such that 
\[
{}^t X A X=           \begin{pmatrix} 
                             \tilde J_r & O\\
                              O      & O
                      \end{pmatrix}.
\]
\end{Lemma}

\proof Induction on $n$.
If $A=O$, then there is nothing to be proved.
Assume that $A\neq O$.
Since $A=(a_{ij})$ 
is alternating, there exists some $1\leq i< j\leq n$ such that 
$a_{ij}\neq 0$.
Exchanging the first row and the $i$th row and then exchanging the first 
column and the $i$th column, we may assume that $a_{1j}\neq 0$.
Exchanging the $j$th column and the $n$th column, and then exchanging
the $j$th row and the $n$th row, we may assume that $a_{1n}\neq 0$.
Multiplying appropriate $X$ from the right and changing $A$ except for
the first column, we may assume that 
$a_{1n}=1$, and $a_{1j}=0$ for $j<n$, and then multiplying ${}^tX$ from the
left, 
$A$ is still alternating.
Multiplying appropriate $Y$ from the left and changing $A$ except for
the first row, we may further assume that $a_{in}=0$ for $i>1$.
Multiplying ${}^t Y$ from the right, $A$ is still alternating.
Thus we may assume that $A$ is of the form
\[
A = \begin{pmatrix}
         0 & 0 & 1\\
         0 & A_1 & 0\\
         -1& 0   &0
       \end{pmatrix},
\]
where $A_1$ is an $(n-2)\times (n-2)$ alternating matrix.
By induction assumption, the rest is easy.
\qed

\begin{Lemma}\label{surjective.thm}
$\pi\colon M\rightarrow Y_t$ is surjective.
\end{Lemma}

\proof Let $\psi\in Y_t$.
Let $\rank\psi=2l$, where $l\leq t-1$.
It is easy to find $Z\in\Mat(2t-2,k)$ such that
\[
{}^t Z \tilde J_{t-1} Z = \begin{pmatrix}
                               \tilde J_l & O\\
                                  O       & O
                                \end{pmatrix}.
\]
By Lemma~\ref{alternating.thm}, we may write
\[
{}^t X A(\psi) X=\begin{pmatrix}
                    \tilde J_{l} & O\\
                        O          & O
                      \end{pmatrix}
\]
with $X\in\GL_n(k)$.
Let $Y$ be the $2l\times n$ matrix consisting of the first $2l$ rows
of $X^{-1}$.
Then we have ${}^t Y \tilde J_{l} Y=A(\psi)$.
If $T=ZY'$, where $Y'$ is the $(2t-2)\times n$ matrix whose first $2l$ rows
are $Y$, and the rest is zero, then ${}^tT \tilde J_{t-1} T=A(\psi)$.
If $\varphi\in\Hom(V,E)$ whose representation matrix is $T$, then
$\pi(\varphi)=\psi$.
\qed

\begin{Lemma}\label{cmint.thm}
$Y_t$ is Cohen--Macaulay and integral.
\end{Lemma}

\proof It is well-known that $k[Y_t]$ is a graded 
ASL on a distributive lattice,
see \cite[section~12]{DEP} and \cite{Kurano}.
Hence $k[Y_t]$ is Cohen--Macaulay by \cite[Theorem~8.1]{DEP}.
It is also reduced by \cite[(5.7)]{BV}.
Since $M$ is irreducible and $\pi$ is surjective, $Y_t$ is irreducible.
Hence $Y_t$ is integral.
\qed

Set $U:= Y_t\setminus Y_{t-1}=\{\varphi\in\Alt(V)\mid \rank \varphi=2(t-1)\}$.

\begin{Lemma} \label{nonsing.thm}
$U$ is non-singular.
\end{Lemma}

\proof By Lemma~\ref{cmint.thm}, $U$ is integral.
On the other hand, $\GL(V)$ acts on $U$ by $(g,u)\mapsto g^* u g$.
By Lemma~\ref{alternating.thm}, this action is transitive.
Hence $U$ is non-singular.
\qed

\begin{Proposition}\label{principal.thm}
The morphism $\pi|_{\pi^{-1}(U)}\colon \pi^{-1}(U)\rightarrow U$ is
a principal $G$-bundle.
\end{Proposition}

\proof Note that $\pi$ is $G$-invariant.
Consider that $U\subset\Alt(V)$.
Note that $\pi^{-1}(U)=\{\varphi\in\Hom(V,E)\mid \varphi 
\text{ surjective}\}$.

Consider the morphism $\Phi\colon G\times \pi^{-1}(U)\rightarrow
\pi^{-1}(U)\times_U \pi^{-1}(U)$.
Let $(\varphi,\psi)\in \pi^{-1}(U)\times_U \pi^{-1}(U)$.
Then $\varphi^*\tilde \rho\varphi=\psi^*\tilde\rho\psi$.
Hence
\[
\Ker\varphi=\Ker(\varphi^*\tilde\rho\varphi)=
\Ker(\psi^*\tilde\rho\psi)=\Ker\psi.
\]
By the homomorphism theorem, there exists some $g\in\GL(E)$ such that
$\varphi=g\psi$.
Since $\psi^*g^*\tilde\rho g \psi=\psi^*\tilde\rho\psi$, $\psi$ is
surjective, and $\psi^*$ is injective, we have that $g^*\tilde\rho g=\tilde
\rho$.
This shows that $g\in G$.
Since $(\varphi,\psi)=\Phi(g,\psi)$, $\Phi$ is surjective.
The injectivity of $\Phi$ is easy.
So $\Phi$ is an isomorphism.

Since $\pi$ is surjective by Lemma~\ref{surjective.thm}, $\pi|_{\pi^{-1}(U)}$
is also surjective.
Since $U$ and $\pi^{-1}(U)$ are non-singular by Lemma~\ref{nonsing.thm},
$\pi|_{\pi^{-1}(U)}$ is faithfully flat as in the proof of 
Theorem~\ref{dp.thm}.
\qed

\begin{Corollary}\label{dimY_t.thm}
$\dim Y_t=(2n-2t+1)(t-1)$.
\end{Corollary}

\proof By the proposition,
\begin{multline*}
\dim Y_t=\dim U=\dim \pi^{-1}(U)-\dim G\\
=2n(t-1)-(2(t-1)^2+(t-1))=
(2n-2t+1)(t-1).
\qed
\end{multline*}

\begin{Corollary}\label{R_4.thm}
$\codim_{Y_t}Y_{t-1}\geq 5$.
\end{Corollary}

\proof By Corollary~\ref{dimY_t.thm}, 
\[
\codim_{Y_t}Y_{t-1}=2(n-2t)+5\geq 5.
\qed
\]

\begin{Corollary}\label{normal.thm}
$Y_t$ is normal.
\end{Corollary}

\proof By Corollary~\ref{R_4.thm} and Lemma~\ref{nonsing.thm}, 
we have that $Y_t$ satisfies the $(R_4)$ condition.
Since $Y_t$ is Cohen--Macaulay by Lemma~\ref{cmint.thm}, 
$Y_t$ is normal.
\qed

\def\citinfo{\cite[(6.6), (6.7)]{DP}}
\begin{Theorem}[\citinfo]
Let $R$ be a commutative ring, and consider the morphism
$\pi\colon M_R \rightarrow (Y_t)_R$ over $R$.
Then the associated ring homomorphism $R[(Y_t)_R]\rightarrow R[M_R]^{G_R}$ 
is an isomorphism.
\end{Theorem}

\proof As a $G$-module, $\Hom(V,E)\cong E^{\oplus n}$.
Note that $E\cong E^*$ as a $G$-module.

By \cite[section~4]{AJ} and \cite[Corollary~III.4.1.8]{Hashimoto}, 
$\Sym E$ has a good filtration as a $G$-module.
By Mathieu's theorem \cite{Mathieu}, $\Bbb Z[M_{\Bbb Z}]$ has a good
filtration.
As in section~\ref{arbitrary.sec}, we may assume that the base ring
$R$ is an algebraically closed field $k$.

Since $Y_t$ is integral and $\pi$ is surjective, the associated map
$k[Y_t]\rightarrow k[M]^G$ is injective.
We want to prove that this is surjective.
Since $Y_t$ is normal by Corollary~\ref{normal.thm} and $\codim_{Y_t}Y_{t-1}
\geq 2$ by Corollary~\ref{R_4.thm}, we have that
$k[Y_t]\cong \Gamma(U,\Cal O_U)$, where $U=Y_t\setminus U_{t-1}$.
Since $\pi|_{\pi^{-1}(U)}\colon \pi^{-1}(U)\rightarrow U$ is a principal
$G$-bundle by Proposition~\ref{principal.thm}, 
$\Gamma(U,\Cal O_U)\rightarrow\Gamma(\pi^{-1}(U),\Cal O_M)^G$ is surjective
by Lemma~\ref{torsor.thm}.
Hence $k[Y_t]\rightarrow k[M]^G$ is surjective, as desired.
\qed

\end{document}